\documentclass{elsarticle}
\usepackage{amsmath}
\usepackage{amssymb}
\usepackage{amsthm}
\usepackage{placeins}
\usepackage{enumitem}
\usepackage{cleveref}
\newtheorem{Theorem}{Theorem}
\newtheorem{Lemma}[Theorem]{Lemma}
\newtheorem{Proposition}[Theorem]{Proposition}
\newtheorem{Corollary}[Theorem]{Corollary}
\newtheorem{Conjecture}[Theorem]{Conjecture}
\newtheorem{Definition}[Theorem]{Definition}
\newtheorem{Remark}[Theorem]{Remark}

\begin{document}

\begin{frontmatter}

\title{On the Termination of the General XL Algorithm and Ordinary Multinomials}

\author{Gary McGuire\fnref{GM}}
\author{Daniela O'Hara\fnref{DO}}
\fntext[GM]{Research supported by Science Foundation Ireland Grant 13/IA/1914.}
\fntext[DO]{Research supported by a Postgraduate Government of Ireland Scholarship from the Irish Research Council.}
\address{University College Dublin, Belfield, Dublin 4, Ireland}

\begin{abstract}
The XL algorithm is an algorithm for solving overdetermined systems of multivariate polynomial  equations, which was initially introduced for quadratic equations. However, the algorithm works for polynomials of any degree, and in this paper we will focus on the performance of XL for 
polynomials of degree $\geq3$, where the optimal termination value of 
the parameter $D$ is still unknown. We prove that the XL algorithm terminates at a certain value of $D$ in the case that the number of equations exceeds the number of variables by 1 or 2.
We also give strong evidence that this
value is best possible, and we show that this value is smaller than the degree of regularity.
Part of our analysis requires proving that ordinary multinomials are strongly unimodal, and this result may be of independent interest.
\end{abstract}


\end{frontmatter}

\section{Introduction}

There are many algorithms and methods for solving systems of polynomial equations in several variables, and the XL algorithm is one such method. It was first used as a tool for attacking the HFE (Hidden Field Equations) cryptosystem. The XL algorithm is part of a family of algorithms, which is based on an idea by Lazard (\cite{Lazard}) that uses linear algebra to solve the system. Other algorithms in this family include F4 and F5 (\cite{F4,F5}). %In fact, the XL algorithm can be seen as a redundant variant of F4 (see \cite{XLF4}).

Multivariate cryptosystems received a lot of attention during the 1990s, and are among the proposals today for post-quantum cryptography.
In multivariate cryptography, one needs to solve a system of multivariate polynomial equations over a finite field in order to recover the original message. Moreover, multivariate polynomial equations have appeared in the cryptanalysis of other public key cryptosystems, such as the Elliptic Curve Discrete Logarithm Problem, where the discovery of Semaev's Summation Polynomials (\cite{Semaev04}) has led to an increased interest is solving such systems.
Thus, cryptography provides one application area for algorithms such as the XL algorithm. 

Most papers so far that analyse the XL algorithm have considered quadratic equations only. In this article we will present some results on the success parameter of the XL algorithm for arbitrary degree polynomials.

In \Cref{xlback} we give background on the XL algorithm, and in \Cref{hilbert} we prove a new theoretical result in commutative algebra which will 
allow us to estimate the optimal termination value of the parameter $D$.

In \Cref{om} we will prove that the ordinary multinomials $\binom{N}{k}_s$ are strongly unimodal, as well as finding the smallest $k$ such that the inequality $\binom{N}{k}_s\leq k$ holds. These results are of independent interest and this section may be read independently of the rest of the paper. We will prove another inequality involving ordinary multinomials in \Cref{c2}. Other papers have outlined proofs of unimodality of ordinary multinomials before. Our method is different and proves \emph{strong} unimodality.

\Cref{c1} contains the main results of the paper -- our theoretical and computational results on the XL algorithm in the case of one more equation than unknown.  \Cref{c2} considers the case of two more equations than unknowns.

One important  concept in the  analysis of algorithms for computing Groebner bases is the degree of regularity,
In Section 7 we will compare the parameter $D$ in this paper 
 to the degree of regularity.
Our investigations show that the XL algorithm succeeds (sometimes)  at a $D$ which is smaller  than the degree of regularity.

\section{Background on the XL algorithm}\label{xlback}

The XL (eXtended Linearization) Algorithm, introduced in \cite{CKPS00}, is an algorithm to solve (overdetermined) systems of multivariate polynomial equations. Consider a system of multivariate polynomial equations
\begin{align*}
 f_1(x_1,\ldots,x_n)&=0\\
 &\ldots\\
 f_{n+c}(x_1,\ldots,x_n)&=0
\end{align*}
over a (finite) field K, where $c\geq1$.

Fix $D\in\mathbb{N}$, where $D>\deg(f_i)$ for all $i$. 
We call $D$ the maximal degree. We consider the system of all products
\begin{equation*}
\prod^k_{l=1} x_{j_l}\cdot f_i(x_1,\ldots,x_n)
\end{equation*}
where $k\leq D-\deg(f_i)$ for $i=1,\ldots,n+c$. Note that \cite{CKPS00} assumes $\deg(f_i)=2$ for all $i$, an assumption we do not make. The idea of the XL algorithm is to linearize this new system in hope of finding a univariate equation which then allows us to solve the initial system.

%Previous versions of the XL algorithm have fixed $D$ beforehand, 
%however we will slightly modify the algorithm to include $D$ as a variable.

\begin{Definition}{(The XL Algorithm)}

Input: polynomials $f_1,\ldots,f_{n+c}$ in variables $x_1,\ldots,x_n$ (an overdetermined system with 0-dimensional solution space) and a positive integer $D$.

Output: All solutions of the system $f_1=0,\ldots,f_{n+c}=0$.

\begin{enumerate}
 \item \textbf{Multiply}: Generate all the products $\prod^k_{l=1}x_{j_l}\cdot f_i(x_1,\ldots,x_n)$ 
 with $k\leq D-deg(f_i)$.
 \item \textbf{Linearize}: Consider each monomial in the $x_i$ of degree $\leq D$ as a new variable 
 and perform Gaussian elimination on the equations obtained in step 1.\\
 The ordering on the monomials must be such that all the terms containing one (fixed) variable (say $x_1$) are eliminated last.
 \item \textbf{Solve}: If step 2 yields at least one univariate equation in the powers of $x_1$, solve this equation over $K$. If not, algorithm fails.
 \item \textbf{Repeat}: Substitute each solution for $x_1$ into the other equations and repeat the process to find the values of the other variables.
\end{enumerate}
\end{Definition}

See \cite{CY04} for some variations and discussions of the algorithm, and a comparison with Gr\"{o}bner basis algorithms. In particular, it was shown by Moh \cite{Moh} (see also \cite{CY04}) that the algorithm terminates for some $D$ provided the solution set
is 0-dimensional.

One might set the starting value of $D$ at $D=1+\max_i \deg(f_i)$. If the ``Solve" step fails for this $D$, or any $D$, we increment the input value of $D$ and run the algorithm again. In practice, the first $D$ at which the algorithm succeeds and terminates is usually larger than $D=1+\max_i\deg(f_i)$.

Considering only those $D$ for which the XL algorithm succeeds, the algorithm increases in running time with $D$. Therefore, when looking for maximum efficiency we would like to know the smallest $D$ such that the XL algorithm succeeds. Let us denote this value by $D^*$. If we know $D^*$, we could use this as the starting input value for $D$.
%If we don't know the exact value of $D^*$ then a lower bound for $D^*$  is useful because we can use the lower bound as the starting value in the XL algorithm. Using ordinary multinomials, we will derive some lower bounds in this article.

\begin{Remark}
Step 2 performs Gaussian elimination. For this, we build a matrix whose columns are indexed by monomials of degree up to D (in the chosen order), and whose rows are indexed by the polynomials generated in step 1. The entries of the matrix in row $f$ 
are the coefficients of the monomials in $f$. This matrix is usually called the Macaulay matrix.
\end{Remark}

\section{A Hilbert series associated to the XL algorithm}\label{hilbert}

We will now follow C. Diem (\cite{Diem04}) to set up the theoretical background for determining the optimal choice of the maximal degree $D$.

Let $V_D$ be the $K$-vector space generated by the products produced in the first step of the XL algorithm, i.e.
\begin{equation}\label{eq:V_D}
V_D=\langle\prod^k_{l=1}x_{j_l}\cdot f_i(x_1,\ldots,x_n)\text{ with }k\leq D-deg(f_i)\rangle_K.
\end{equation}
Let $K[x_1,\ldots,x_n]_{\leq D}$ be the $K$-vector space of polynomials of total degree $\leq D$. Let
\begin{equation}\label{eq:chi(D)}
\chi(D):=dim_K(K[x_1,\ldots,x_n]_{\leq D})-dim_K(V_D).
\end{equation}

\begin{Theorem}\label{success} (\cite{Diem04})
If $\chi(D)\leq D$ then the XL algorithm terminates for that $D$.
\end{Theorem}

\begin{proof}
Since $dim_K(K[x_1]_{\leq D})=D+1$, if $\chi(D)\leq D$ then $dim_K(K[x_1]_{\leq D})>\chi(D)$, i.e. $dim_K(V_D)+dim_K(K[x_1]_{\leq D})>dim_K(K[x_1,\ldots,x_n]_{\leq D})$ and hence $V_D\cap K[x_1]_{\leq D}\neq \{0\}$. So step 2 of the XL algorithm produces a univariate equation in $x_1$ in this case.
\end{proof}

Clearly then, for reasons of efficiency of the XL algorithm, we would like to know the the smallest $D$ such that $\chi(D)\leq D$. We will now try to estimate $\chi(D)$ in order to find the smallest $D$ such that $\chi(D)\leq D$.

Let $K[x_0,x_1,\ldots,x_n]_D$ be the $K$-vector space of all homogeneous polynomials of total degree $D$. Let $F_i\in K[x_0,x_1,\ldots,x_n]$ denote the homogenization of $f_i$. Then $V_D\cong I_D$ via the degree D homogenization map, where
\begin{equation*}
 I_D := \langle \prod^k_{l=1} x_{j_l}\cdot F_i(x_0,x_1,\ldots,x_n) \text{ with } k=D-deg(F_i) \rangle_K.
\end{equation*}
Now $I_D$ is the $D^{\text{th}}$ homogeneous component of the homogeneous ideal $I:=(F_1,\ldots,F_{n+c})\vartriangleleft K[x_0,\ldots,x_n]$. It follows that
\begin{align*}
 \chi(D)&=dim_K(K[x_0,\ldots,x_n]_D)-dim_K(I_D)\\
 &=dim_K(K[x_0,\ldots,x_n]_D/I_D)\\
 &=dim_K((K[x_0,\ldots,x_n]/I)_D).
\end{align*}

\begin{Definition}
Let $R:=K[x_0,\ldots,x_n]$. Let $M=\oplus_{i\in\mathbb{N}}M_i$ be a finitely generated positively
graded $R$ module. Define the Hilbert function of $M$ as
\begin{align*}
\chi_M :\text{ } \mathbb{N}_0&\rightarrow\mathbb{N}_0\\
i&\mapsto dim_R(M_i)
\end{align*}
and the Hilbert series of $M$ as $H_M:=\displaystyle\sum_{i\in\mathbb{N}}\chi_M(i)T^i$.
\end{Definition}

So $\chi(D)=\chi_{R/I}(D)$ $\forall D\in\mathbb{N}$.

\begin{Definition}{(see also \cite{Froberg85})}
A form $G\in K[x_0,\ldots,x_n]$ of degree $d$ is called \textit{generic} if all monomials of degree $d$ in $K[x_0,\ldots,x_n]$ have non-zero coefficients in $G$, and those coefficients are algebraically independent over the prime field of $K$.
\end{Definition}

\begin{Lemma}
The Hilbert series of an ideal generated by a generic system of forms depends only on the field characteristic, the rank $n$, and the degrees of the forms. If $\text{char } K=0$, we call it the generic Hilbert series of type $(n+1;m;d_1,\ldots,d_m)$ where
$m$ is the number of forms. (For us, $m=n+c$.)
\end{Lemma}

\begin{proof}
See \cite{Diem04}.
\end{proof}

\begin{Proposition}\label{Hilbert_series_type}
Let $K$ be a field (any characteristic), let $F_1,\ldots,F_m\in R:=K[x_0,\ldots,x_n]$ be forms of degree $d_1,\ldots,d_m$ (not necessarily generic). Let $I:=(F_1,\ldots,F_m)\vartriangleleft R$. Let $H_g$ be the generic Hilbert series of type $(n+1;m;d_1,\ldots,d_m)$. Then $H_{R/I}\geq H_g$ coefficient-wise.
\end{Proposition}

\begin{proof}
See \cite{Diem04}.
\end{proof}

\begin{Proposition}\label{multiplicationmap}
Let $G_1,\ldots,G_m=G\in R$ be a generic system of forms with $m\leq n+1$, let $d=deg(G)$. Let $J:=(G_1,\ldots,G_{m-1})\vartriangleleft R$. Then for all $D\in\mathbb{N}_0$, the multiplication map
\begin{align*}
G\cdot : (R/J)_D&\rightarrow(R/J)_{D+d}\\
\bar{F}&\mapsto G\cdot \bar{F}
\end{align*}
is injective, and we have a short exact sequence
\begin{align*}
0 \longrightarrow (R/J)_D \xrightarrow{G\cdot} (R/J)_{D+d} \longrightarrow (R/(J,G))_{D+d}\longrightarrow 0.
\end{align*}
\end{Proposition}

\begin{Proposition}\label{generic m<=n+1}
Let $G_1,\ldots,G_m\in R$ be a generic system of forms of degrees $d_1,\ldots,d_m$ with $m\leq n+1$. Then 
\begin{align*}
H_{R/(G_1,\ldots,G_m)}=\frac{\prod_{j=1}^m(1-T^{d_j})}{(1-T)^{n+1}} \longrightarrow 0.
\end{align*}
\end{Proposition}

\begin{proof}
By the previous proposition, $H_{(R/J)_{D+d}}=H_{(R/J)_D}+H_{(R/(J,G))_{D+d}}$. The result then follows by induction since $H_R=\frac{1}{(1-T)^{n+1}}$.
\end{proof}

All preceding results of this section can be found in \cite{Diem04}. We incude the statements because they are needed for the next result, which is new.

\begin{Theorem}\label{Hilbert_inequality}
Let $F_1,\ldots,F_{n+c}\in R$ $(\text{with }c\geq1)$ be forms of degrees $d_1,\ldots,d_{n+c}$. Then
\begin{align*}
H_{R/(F_1,\ldots,F_{n+c})}\geq\Big(1-\displaystyle\sum_{j=n+2}^{n+c}T^{d_j}\Big)
\frac{\prod_{j=1}^{n+1}(1-T^{d_j})}{(1-T)^{n+1}}
\end{align*}
 coefficient-wise. If $d_1=\ldots=d_{n+c}=d$ then 
\begin{align*}
H_{R/(F_1,\ldots,F_{n+c})}\geq(1-(c-1)T^d)(1+T+\ldots+T^{d-1})^{n+1}
\end{align*}
coefficient-wise.
\end{Theorem}

\begin{proof}
The case $d_1=\ldots=d_{n+c}=2$ is done in \cite{Diem04}. The general case is similar. We will prove that the generic Hilbert series of type $(n+1;n+c;d_1,\ldots,d_{n+c})$ is $\geq\Big(1-\displaystyle\sum_{j=n+2}^{n+c}T^{d_j}\Big)\frac{\prod_{j=1}^{n+1}(1-T^{d_j})}{(1-T)^{n+1}}$ coefficient-wise. The result then follows by \Cref{Hilbert_series_type}.

Let $G_1,\ldots,G_{n+c}\in R:=K[x_0,\ldots,x_n]$ (with $\text{char } K=0$) be a generic system of forms of degrees $d_1,\ldots,d_{n+c}$. Let $R':=R/(G_1,\ldots,G_{n+1})$ and let $I':=(G_{n+2},\ldots,G_{n+c})\vartriangleleft R'$. Then $H_{R'}=\frac{\prod_{j=1}^{n+1}(1-T^{d_j})}{(1-T)^{n+1}}$ by \Cref{generic m<=n+1}. 
Since $R/(G_1,\ldots,G_{n+c})\cong R'/I'$, 
\begin{align*}
\chi_{R/(G_1,\ldots,G_{n+c})}(D)=\chi_{R'/I'}(D)=dim_K(R'_D)-dim_K(I'_D).
\end{align*}
For $D\geq d$, $I'_D=\sum_{j=n+2}^{n+c}G_j\cdot R'_{D-d_j}$ with multiplication map $G_j\cdot: R'_{D-d_j}\rightarrow R'_D$ as in \Cref{multiplicationmap}. Hence 
$dim_K(I'_D)\leq\sum_{j=n+2}^{n+c}dim_K(R'_{D-d_j})$ by \Cref{multiplicationmap}. Thus we have 
\begin{align*}
\chi_{R/(G_1,\ldots,G_{n+c})}(D)\geq dim_K(R'_D)-\displaystyle\sum_{j=n+2}^{n+c}dim_K(R'_{D-d_j})
\end{align*}
and hence
\begin{align*}
H_{R/(G_1,\ldots,G_{n+c})}\geq H_{R'}-\displaystyle\sum_{j=n+2}^{n+c}T^{d_j}H_{R'}=
\Big(1-\displaystyle\sum_{j=n+2}^{n+c}T^{d_j}\Big)\frac{\prod_{j=1}^{n+1}(1-T^{d_j})}{(1-T)^{n+1}}.
\end{align*}
\end{proof}

\begin{Corollary}\label{xi_inequality}
If $f_1,\ldots,f_{n+1}\in K[x_1,\ldots,x_n]$ (i.e. $c=1$) with $deg(f_i)=d$ for all $i$, then $\chi(D)\geq \binom{n+1}{D}_{d-1}$. Furthermore, if all the polynomials are generic, then equality holds.
%with $\chi(D)$ as defined in (\cref{eq:chi(D)}) at the start of this section.
\end{Corollary}

\begin{proof}
Take $c=1$ in \Cref{Hilbert_inequality}, and recall that the coefficient of $T^D$ in $(1+T+\ldots+T^{d-1})^{n+1}$ is $\binom{n+1}{D}_{d-1}$. The second statement follows from \Cref{generic m<=n+1}.
\end{proof}

In \Cref{c1} we use this Corollary to find the smallest $D$ such that $\chi(D)\leq D$ (and XL succeeds
by \Cref{success}) in the $c=1$ case.

\begin{Remark}\label{Ds}
We point out that there are three possibly different $D$'s under discussion:
\begin{enumerate}
 \item the smallest $D$ such that the XL algorithm terminates, denoted by $D^*$,
 \item the smallest $D$ such that $\chi(D) \leq D$, denoted by $D_\chi$,
 \item the smallest $D$ such that $\binom{n+1}{D}_{d-1}\leq D$, denoted by $D_m$.
\end{enumerate}
\noindent \Cref{success} implies that $D^*\leq D_\chi$. 
In the case $c=1$ and all degrees equal, 
\Cref{xi_inequality} implies that $D_m\leq D_\chi$, and $D_m= D_\chi$ when the polynomials are generic. %So for generic polynomials, we have $D^*\leq D_\chi=D_m$.

In \Cref{c1} we investigate the relationship between these $D$s when $c=1$. At the end of the section we will conjecture that $D^*=D_\chi=D_m$ when $c=1$ and the polynomials are generic, and provide evidence.
\end{Remark}

\begin{Remark}\label{Ds2}
It is intuitively obvious that generic polynomials lead to the worst-case $D^*$, i.e., 
the XL algorithm with non-generic polynomials will
terminate with a smaller $D$ than for generic polynomials. 
Assuming this to be true, for non-generic polynomials we have $D^*\leq D_m\leq D_\chi$ when $c=1$
and all degrees are equal. This agrees with our experiments, see \Cref{tab:exp}.
\end{Remark}

\section{Ordinary Multinomials}\label{om}

This section is independent of the rest of the paper. 
We will prove here that ordinary multinomials are strongly unimodal, and some inequalities.

\begin{Definition}
A sequence $s_0, s_1, \ldots ,s_N$ of integers is said to be \emph{unimodal} if there is an integer $t$ with $0\leq t \leq N$ such that
\begin{align*}
s_0\leq s_1\leq \cdots \leq s_t, \quad s_t\geq s_{t+1}\geq \cdots \geq s_N.
\end{align*}
\end{Definition}

A unimodal sequence is said to be \emph{strongly unimodal} if all the inequalities are strict, except that $s_t=s_{t+1}$ may hold. For example, the sequence of binomial coefficients $\binom{N}{k}$ ($k=0,1,\ldots,N$) is strongly unimodal.

\begin{Remark}
We remark that there are different definitions of \emph{strongly unimodal} in the literature. One definition is the same as ours except it does not allow $s_t=s_{t+1}$, in which case the binomial coefficients are strongly unimodal only for $N$ even.
\end{Remark}

\bigskip

\begin{Definition}
The coefficient of $T^k$ in $(1+T+\ldots+T^{s})^{N}$ is called \textit{ordinary multinomial} or \textit{generalized binomial coefficient of order $s$} and is denoted by $\binom{N}{k}_s$.
\end{Definition}

Combinatorially, $\binom{N}{k}_s$ is the number of different ways of distributing $k$ objects among $N$ 
boxes, where each box contains at most $s$ objects.
See \cite{Boll86} and \cite{Bondarenko} for an introduction.

\begin{table}[ht!]
\setlength{\tabcolsep}{2pt}
\begin{tabular}{rccccccccccccccccccccccccc}
$N=0$:& & & & &  & &  & & &  & &  &1\\
$N=1$:& & & & &  & &  & & &1 & &1 & &1 & &1\\
$N=2$:& & & & &  & &1 & &2&  &3&  &4&  &3&  &2& &1\\
$N=3$:& & & &1&  &3&  &6& &10& &12& &12& &10& &6&  &3&  &1\\
$N=4$:&1& &4& &10& &20& &31& &40& &44& &40& &31&&20& &10& &4& &1\\
\end{tabular}
\caption{\label{tab:expl}Triangle of ordinary multinomials $\binom{N}{k}_s$ with $s=3$ and $k=0,...,sN$.}
\end{table}

\begin{Proposition}\label{ordinarymultinomials}
\begin{enumerate}[label=(\alph*),ref=(\alph*)]
 \item $\binom{N}{k}_s=0$ for $k<0$ and $k>sN$.
 \item For $s=1$, the ordinary multinomials are the usual binomial coefficients.
 \item $\dbinom{N}{k}_s=\displaystyle\sum_{i=0}^{\lfloor k/(s+1)\rfloor}(-1)^i\binom{N}{i}
 \binom{k-i(s+1)+N-1}{k-i(s+1)}$
 \item $\dbinom{N}{k}_s=\displaystyle\sum_{\substack{j_1+\ldots+j_{s+1}=N\\
   j_2+2j_3+\ldots+sj_{s+1}=k}}\binom{N}{j_1,\ldots,j_{s+1}}$.
\end{enumerate}
\end{Proposition}

\begin{proof}
(c) follows from expanding $(1+T+\ldots+T^{s})^{N}=\frac{(1-T^{s+1})^N}{(1-T)^{N}}$.\\
(d) follows from the multinomial theorem.
\end{proof}

%The XL algorithm succeeds when $\chi(D)\leq D$, so we will try to find the smallest $D$ such that $\binom{n+1}{D}_{d-1}\leq D$ since this gives a lower bound for the case $c=1$ by \cref{xi_inequality}.

\begin{Lemma}
$\binom{N}{k}_s=\binom{N}{sN-k}_s$, i.e. the ordinary multinomials are symmetric.
\end{Lemma}

\begin{proof}
Using \Cref{ordinarymultinomials}(d), if $j_1,\ldots,j_{s+1}$ satisfy $\{\substack{j_1+\ldots+j_{s+1}=N\\j_2+2j_3+\ldots+sj_{s+1}=k}$, then $j_s+2j_{s-1}+\ldots+sj_{1}=sN-k$ and $\binom{N}{j_{s+1},\ldots,j_1}=\binom{N}{j_1,\ldots,j_{s+1}}$.
\end{proof}

\begin{Lemma}\label{recurrence}
$\binom{N}{k}_s=\sum_{m=0}^s \binom{N-1}{k-m}_s$.
\end{Lemma}

\begin{proof}
Use \Cref{ordinarymultinomials}(d) and the recurrence relation of multinomial coefficients.
\end{proof}

\begin{Theorem}\label{strongunimodality}
The ordinary multinomials are strongly unimodal. For $N\geq2$ we have
\begin{align*}
\binom{N}{0}_s<\binom{N}{1}_s<\cdots<\binom{N}{\lfloor\frac{sN}{2}\rfloor}_s=
\binom{N}{\lceil\frac{sN}{2}\rceil}_s>\cdots>\binom{N}{sN}_s.
\end{align*}
\end{Theorem}

\begin{proof}
It follows from \Cref{recurrence} that
\begin{align*}
\dbinom{N}{k+1}_s&=\displaystyle\sum_{m=0}^s\binom{N-1}{k+1-m}_s\\
&=\dbinom{N-1}{k+1}_s+\displaystyle\sum_{m=0}^s\binom{N-1}{k-m}_s-\dbinom{N-1}{k-s}_s\\
&=\dbinom{N-1}{k+1}_s+\dbinom{N}{k}_s-\dbinom{N-1}{k-s}_s.
\end{align*}
So if $\binom{N-1}{k+1}_s>\binom{N-1}{k-s}_s$ for $k<\lfloor\frac{sN}{2}\rfloor$ then $\binom{N}{k+1}_s>\binom{N}{k}_s$ for $k<\lfloor\frac{sN}{2}\rfloor$. We proceed by induction on $N$.

For $N=2$, the base case, $\binom{1}{k+1}_s-\binom{1}{k-s}_s=1-0=1$ for $0\leq k<s$. Thus $\binom{2}{k+1}_s>\binom{2}{k}_s$ for all $k<\lfloor\frac{sN}{2}\rfloor$.

Now assume $\binom{N-1}{k+1}_s>\binom{N-1}{k}_s$ for all $k<\lfloor\frac{s(N-1)}{2}\rfloor$. Then 
\begin{align*}
\binom{N-1}{k+1}_s>\binom{N-1}{k}_s>\cdots>\binom{N-1}{k-s}_s
\end{align*}
whenever $k<\lfloor\frac{s(N-1)}{2}\rfloor$.
It remains to prove the cases $\lfloor\frac{s(N-1)}{2}\rfloor \leq k <\lfloor\frac{sN}{2}\rfloor$.

If $k=\lfloor\frac{s(N-1)}{2}\rfloor+m$ with $0\leq m<s/2$ (and hence $k<\lfloor\frac{sN}{2}\rfloor$) then
\begin{align*}
\binom{N-1}{k+1}_s=\binom{N-1}{s(N-1)-k-1}_s=\binom{N-1}{\lceil\frac{s(N-1)}{2}\rceil-m-1}_s.
\end{align*}

\underline{Case 1:} where $ s(N-1)$ is even. In this case $\lceil\frac{s(N-1)}{2}\rceil=\lfloor\frac{s(N-1)}{2}\rfloor$ so
\begin{align*}
\binom{N-1}{k+1}_s&=\binom{N-1}{\lfloor\frac{s(N-1)}{2}\rfloor-m-1}_s>\cdots\\
&>\binom{N-1}{\lfloor\frac{s(N-1)}{2}\rfloor+m-s}_s=\binom{N-1}{k-s}_s
\end{align*}
by induction assumption and since $m<s/2$.

\underline{Case 2:} where $ s(N-1)$ is odd. In this case $\lceil\frac{s(N-1)}{2}\rceil=\lfloor\frac{s(N-1)}{2}\rfloor+1$ so
\begin{align*}
\binom{N-1}{k+1}_s&=\binom{N-1}{\lfloor\frac{s(N-1)}{2}\rfloor+1-m-1}_s
=\binom{N-1}{\lfloor\frac{s(N-1)}{2}\rfloor-m}_s\\
&>\cdots>\binom{N-1}{\lfloor\frac{s(N-1)}{2}\rfloor+m-s}_s=\binom{N-1}{k-s}_s
\end{align*}
by induction assumption and since $m<s/2$.

Thus, $\binom{N-1}{k+1}_s>\binom{N-1}{k-s}_s$ and hence $\binom{N}{k+1}_s>\binom{N}{k}_s$ for all $k<\lfloor\frac{sN}{2}\rfloor$.

By symmetry, $\binom{N}{k+1}_s<\binom{N}{k}_s$ for all $k\geq\lceil\frac{sN}{2}\rceil$, and $\binom{N}{\lfloor\frac{sN}{2}\rfloor}_s=\binom{N}{\lceil\frac{sN}{2}\rceil}_s$.
\end{proof}

\begin{Remark}
The smallest mode is given by $\lfloor\frac{sN}{2}\rfloor$. If $sN$ is odd, then we have a plateau of two 
modes: $\frac{sN-1}{2}$ and $\frac{sN+1}{2}$. If $sN$ is even, then we have a peak at $\frac{sN}{2}$.
\end{Remark}

\begin{Remark}
One can also prove (weak) unimodality of $\binom{N}{k}_s$ using Theorem 4.7 of \cite{conv} which says that 
the convolution of two symmetric discrete unimodal distributions is again unimodal, along with $\binom{N}{k}_s=\sum_{m=0}^N \binom{N}{m}\binom{m}{k-m}_{s-1}$ (this relation can be proved using $\binom{N}{j_1,\ldots,j_m}=\prod_{i=1}^m \binom{\sum_{l=1}^i j_l}{j_i}$ and \Cref{ordinarymultinomials}(d)). See also \cite{uniproof}.
\end{Remark}

\begin{Proposition}\label{starsandbars}
$\binom{N}{k}_s=\binom{N+k-1}{k}$ for $k\leq s$.
\end{Proposition}

\begin{proof}
$\binom{N+k-1}{k}$ counts the number of ways of putting $k$ objects into $N$ boxes. As noted in \cite{Bondarenko}, $\binom{N}{k}_s$ is the number of ways of putting $k$ objects into $N$ boxes, where each box contains at most $s$ objects.
\end{proof}

\begin{Remark}
By symmetry, if $k\geq s(N-1)$ then $\binom{N}{k}_s=\binom{(s+1)N-k-1}{sN-k}$.
\end{Remark}

Here is a theorem we will need later in Section 5.

\begin{Theorem}\label{smallestk}
For $2\leq s< \frac{N+1}{2}+\frac{2}{N}$, the smallest $k$ such that $\binom{N}{k}_s\leq k$ is $k=sN-1$.
\end{Theorem}

\begin{proof}
\underline{Claim:} $\binom{N}{k}_s>k$ for all $k<\lfloor\frac{sN}{2}\rfloor$.\\
\underline{Proof of claim:} $\binom{N}{0}_s=1>0$. Now assume $\binom{N}{k}_s>k$. Then $\binom{N}{k+1}_s\geq\binom{N}{k}_s+1>k+1$ by strong unimodality (\Cref{strongunimodality}) for $k<\lfloor\frac{sN}{2}\rfloor$.\\
\underline{Claim:} $\binom{N}{k}_s>k$ for all $k\leq sN-2$ and $2\leq s< \frac{N+1}{2}+\frac{2}{N}$.\\
\underline{Proof of claim:} If $\binom{N}{k}_s\leq k$ for some $k\geq\lceil\frac{sN}{2}\rceil$, then 
$\binom{N}{k+1}_s\leq \binom{N}{k}_s+1\leq k+1$ by strong unimodality. By \Cref{starsandbars}, $\binom{N}{sN-2}_s=\binom{N+1}{2}>sN-2$ for $2\leq s< \frac{N+1}{2}+\frac{2}{N}$.\\
Again by \Cref{starsandbars}, $\binom{N}{sN-1}_s=\binom{N}{1}=N\leq sN-1$ for $s\geq2$. Thus, $k=sN-1$ is the smallest $k$ such that $\binom{N}{k}_s\leq k$.
\end{proof}

\begin{Remark}
For $s=1$, the smallest $k$ such that $\binom{N}{k}_1\leq k$ is $k=N$. If $s\geq\frac{N+1}{2}+\frac{2}{N}$, then the smallest $k$ is $<sN-1$. These instances can easily be evaluated directly, see \Cref{tab:Dmin}.
\end{Remark}

We will also use the following result, which will be proved in Section 6.

\begin{Theorem}
The smallest $k$ such that $\binom{N}{k}_{s}-\binom{N}{k-(s+1)}_{s}\leq k$ is 
\begin{enumerate}[label=(\alph*),ref=(\alph*)]
 \item $k=s+1$ for $N=2$,
 \item $k=2s$ for $N=3$,
 \item $k=\lfloor\frac{s(N+1)}{2}\rfloor+1$ for $N\geq4$.
\end{enumerate}
\end{Theorem}

\section{One more equation than unknowns (c=1) in XL}\label{c1}

We are now ready to find our lower bound on  $D_\chi$   for the case $c=1$. By \Cref{success} the XL algorithm will succeed for this $D$, so this may be a good choice for the initial value of $D$. We will see that this is actually very often the optimal choice.

According to \Cref{xi_inequality}, $\chi(D)\geq \binom{n+1}{D}_{d-1}$ when we have $n$ unknowns, $n+1$ equations, and all equations  have the same degree $d$; and equality holds if the polynomials are generic.
%(If the equations do not all have the same degree, we can take $d$ to be the maximum of the degrees.)
According to \Cref{success}, $\chi(D)\leq D$ is a sufficient condition for the XL algorithm to succeed.
Combining these we get
\begin{align*}
D\geq \chi(D)\geq \binom{n+1}{D}_{d-1}
\end{align*}
as sufficient for success, and so we consider the inequality 
\begin{align*}
D\geq  \binom{n+1}{D}_{d-1}
\end{align*}
and we will apply  \Cref{smallestk}.

\begin{Corollary}\label{lb}
Suppose
$f_1,\ldots,f_{n+1}\in K[x_1,\ldots,x_n]$ with $deg(f_i)=d$ for all $i$ where $3\leq d\leq\frac{n+2}{2}+\frac{2}{n+1}$. 
Then:
\begin{enumerate}
\item $D_{m}= (d-1)(n+1)-1$.
\item  $D_{\chi}\ge (d-1)(n+1)-1$.
\item If the polynomials are generic, then  $D_{\chi}=(d-1)(n+1)-1$.
\item If the polynomials are generic, then $D^*\le (d-1)(n+1)-1$.
\end{enumerate}
\end{Corollary}

\begin{proof}
For part 1, take $s=d-1$ and $N=n+1$ in \Cref{smallestk}.
Parts 2, 3 and 4 follow from  \Cref{Ds}. 
\end{proof}

%\begin{Remark}
%This is why we proved \Cref{smallestk} earlier. We note that the proof of \Cref{smallestk} uses strong unimodality of ordinary multinomials, and this is why we included that result in our paper.
%\end{Remark}

%\begin{Remark}
%\Cref{lb} tells us that $D_m=(d-1)(n+1)-1$. By \Cref{Ds}, we know that $D^*\leq D_m$. So \Cref{lb} implies $D^*\le (d-1)(n+1)-1$.
%\end{Remark}

\begin{table}[ht!]
\begin{center}
\begin{tabular}{|c|c|c||c|c|c||c|c|c|}
\hline
$d$ & $n$ & $D$ & $d$ & $n$ & $D$ & $d$ & $n$ & $D$\\
\hline
2 & $\geq2$ & $n+1$ & 6 & $\geq8$ & $5n+4$ & 9 & $2$ & $8n+4$\\
3 & $\geq2$ & $2n+1$ & 7 & $2,3$ & $6n+3$ & 9 & $3,4$ & $8n+5$\\
4 & $2,3$ & $3n+1$ & 7 & $4,\ldots,9$ & $6n+4$ & 9 & $5,\ldots,13$ & $8n+6$\\
4 & $\geq4$ & $3n+2$ & 7 & $\geq10$ & $6n+5$ & 9 & $\geq14$ & $8n+7$\\
5 & $2,\ldots,5$ & $4n+2$ & 8 & $2$ & $7n+3$ & 10 & $2$ & $9n+4$\\
5 & $\geq6$ & $4n+3$ & 8 & $3$ & $7n+4$ & 10 & $3,4$ & $9n+6$\\
6 & $2$ & $5n+2$ & 8 & $4,\ldots,11$ & $7n+5$ & 10 & $5,\ldots,15$ & $9n+7$\\
6 & $3,\ldots,7$ & $5n+3$ & 8 & $\geq12$ & $7n+6$ & 10 & $\geq16$ & $9n+8$\\
\hline
\end{tabular}
\caption{\label{tab:Dmin}The smallest $D$ such that $\binom{n+1}{D}_{d-1}\leq D$}
\end{center}
\end{table}

\begin{table}[ht!]
\begin{center}
\begin{tabular}{|c|c|c|c|c||c|c|c|c|c|}
\hline
$p$&$d$&$n$&$D_{average}$&$D_{min}$&$p$&$d$&$n$& $D_{average}$&$D_{min}$\\ 
 \hline
 3109&2&2&3.00&3 & 3109&5&2&10.00&10\\
 5011&2&2&3.00&3 & 5011&5&2&10.00&10\\
 3109&2&3&4.00&4 & 3109&5&3&14.00&14\\
 5011&2&3&4.00&4 & 5011&5&3&14.00&14\\
 3109&2&4&5.00&5 & 3109&5&4&18.00&18\\
 5011&2&4&5.00&5 & 5011&5&4&18.00&18\\
 \cline{6-10}
 3109&2&5&6.00&6 & 3109&6&2&12.00&12\\
 5011&2&5&6.00&6 & 5011&6&2&12.00&12\\
 3109&2&6&6.99&7 & 3109&6&3&18.00&18\\
 5011&2&6&7.00&7 & 5011&6&3&18.00&18\\
 \cline{6-10}
 3109&2&7&8.00&8 & 3109&7&2&15.00&15\\
 5011&2&7&8.00&8 & 5011&7&2&15.00&15\\
 \cline{1-5}
 3109&3&2&5.00&5 & 3109&7&3&21.00&21\\
 5011&3&2&5.00&5 & 5011&7&3&21.00&21\\
 \cline{6-10}
 3109&3&3&7.00&7 & 3109&8&2&17.00&17\\
 5011&3&3&7.00&7 & 5011&8&2&17.00&17\\
 3109&3&4&9.00&9 & 3109&8&3&25.00&25\\
 5011&3&4&9.00&9 & 5011&8&3&25.00&25\\
 \cline{6-10}
 3109&3&5&11.00&11 & 3109&9&2&20.00&20\\
 5011&3&5&11.00&11 & 5011&9&2&20.00&20\\
 \cline{1-5}
 3109&4&2&7.00&7 & 3109&9&3&29.00&29\\
 5011&4&2&7.00&7 & 5011&9&3&29.00&29\\
 \cline{6-10}
 3109&4&3&10.00&10 & 3109&10&2&22.00&22\\
 5011&4&3&10.00&10 & 5011&10&2&22.00&22\\
 3109&4&4&14.00&14 & 3109&10&3&33.00&33\\
 5011&4&4&14.00&14 & 5011&10&3&33.00&33\\
\hline
\end{tabular}
\caption{\label{tab:exp}Running XL on randomly chosen polynomials over prime fields of order $p$ with one 
more equation than unknowns. Here $D_{average}$ is the average value of the smallest $D$ such that the XL algorithm terminates over 100 experiments with the same parameters, and $D_{min}$ has been calculated according to \Cref{tab:Dmin}. All experiments were done 
using Magma V2.21-6.}
\end{center}
\end{table}

\begin{Remark}
We ran experiments for small values of $d$ and $n$, summarised in \Cref{tab:exp}. They show that 
$D_m=D^*$, i.e. the upper bound for $D^*$ was met with equality in all (but one, see next Remark) of our experiments.
\end{Remark}

\begin{Remark}
In one experiment, namely $p=3109$, $d=2$, $n=6$, the reader will notice that $D_{average}=6.99$ and not 7. This means that 1 out of 100 randomly chosen sets of input polynomials finished with $D=6$ and not $D=7$.
We checked this set of input polynomials, and two coefficients were equal in one polynomial, i.e., the polynomials were not generic in this one case. This is an example that shows that our bound is (conjecturally) tight for generic polynomials, but non-generic polynomials may finish with a smaller $D$.
We mentioned this point earlier in \Cref{Ds2}.
\end{Remark}

Based on the evidence in \Cref{tab:exp}, we now conjecture that our bound in \Cref{lb} is tight for generic polynomials.

\begin{Conjecture}\label{conj32}
When $c=1$ and the polynomials are generic of degree $d$, we have $D^*=(d-1)(n+1)-1$.
\end{Conjecture}

\begin{Remark}
In the notation of \Cref{Ds} we are conjecturing that $D^*=D_\chi=D_m$.
\end{Remark}

\begin{Remark}
Our conjecture is consistent with Proposition 6 in \cite{Diem04} where it is shown that $D_\chi \geq n/(1+\sqrt{c-1})$ for quadratic equations. When $c=1$ this bound becomes $D_\chi \geq n$, and our conjecture with $d=2$ becomes $D^*=n$.
\end{Remark}

\begin{Remark}
When the polynomials do not all have the same degree, we can take $d$ to be the maximum of the degrees and the above arguments will give a good starting input value for $D$, but not necessarily the optimal value.
\end{Remark}

\section{Two more equations than unknowns (c=2) in XL}\label{c2}

When $c>1$ ($c$ being the difference between the number of equations and the number of unknowns), then \Cref{Hilbert_inequality} tells us that $$\chi(D)\geq \binom{n+1}{D}_{d-1}-(c-1)\binom{n+1}{D-d}_{d-1}$$ when we have $n$ unknowns, and all equations have the same degree $d$. When $c=2$ we need the following result.

\begin{Theorem}
The smallest $k$ such that $\binom{N}{k}_{s}-\binom{N}{k-(s+1)}_{s}\leq k$ is 
\begin{enumerate}[label=(\alph*),ref=(\alph*)]
 \item $k=s+1$ for $N=2$,
 \item $k=2s$ for $N=3$,
 \item $k=\lfloor\frac{s(N+1)}{2}\rfloor+1$ for $N\geq4$.
\end{enumerate}
\end{Theorem}

\begin{proof}
\begin{enumerate}[label=(\alph*),ref=(\alph*)]
 \item $N=2$: If $k\leq s$ then $\binom{2}{k}_{s}\stackrel{\Cref{starsandbars}}{=}\binom{2+k-1}{k}=k+1$, and thus $\binom{2}{k}_{s}-\binom{2}{k-(s+1)}_{s}=k+1-0>k$, $\forall k\leq s$. Now $\binom{2}{s+1}_{s}\stackrel{\text{sym.}}{=}\binom{2}{s-1}_{s}\stackrel{\Cref{starsandbars}}{=}\binom{2+s-2}{s-1}=s$ since we have a peak at $s$, and thus $\binom{2}{s+1}_{s}-\binom{2}{s+1-(s+1)}_{s}=s-1\leq s+1$.
 \item $N=3$: If $k\leq s$ then $\dbinom{3}{k}_{s}-\dbinom{3}{k-(s+1)}_{s}\stackrel{\Cref{recurrence}}{=}\displaystyle\sum_{m=0}^k\dbinom{2}{k-m}_{s}-0\geq
 \dbinom{2}{k}_{s}>k$ by part (a) and since the ordinary multinomials are nonnegative. Now let $k=s+t$
 where $1\leq t\leq s$, i.e. $s<k\leq 2s$. Thus
 \begin{align*}
  \binom{3}{k}_{s}&-\binom{3}{k-(s+1)}_{s}
  =\binom{3}{s+t}_{s}-\binom{3}{t-1}_{s}\\
  &\stackrel{\Cref{recurrence}}{=}\sum_{m=0}^s\binom{2}{s+t-m}_{s}-\sum_{m=0}^{t-1}\binom{2}{t-1-m}_{s}\\
  &=\sum_{m=0}^{t-1}\binom{2}{s+t-m}_{s}+\sum_{m=t}^s\binom{2}{s+t-m}_{s}-
  \sum_{m=0}^{t-1}\binom{2}{t-1-m}_{s}\\
  &\stackrel{\text{sym.}}{=}\sum_{m=0}^{t-1}\binom{2}{s-t+m}_{s}+\sum_{m=t}^s\binom{2}{s+t-m}_{s}-
  \sum_{m=0}^{t-1}\binom{2}{t-1-m}_{s}\\
  &\stackrel{\Cref{starsandbars}}{=}\sum_{m=0}^{t-1}\binom{s-t+m+1}{s-t+m}+
  \sum_{m=t}^s\binom{s+t-m+1}{s+t-m}-\sum_{m=0}^{t-1}\binom{t-m}{t-1-m}\\
  &=\sum_{m=0}^{t-1}(s-t+m+1)+\sum_{m=t}^s(s+t-m+1)-\sum_{m=0}^{t-1}(t-m)\\
  &=(s-t+1)(s+2t+1)+\frac{t(t-3)}{2}-\frac{s(s+1)}{2}.
 \end{align*}
 Now we are left with considering the inequality
 \begin{align*}
  (s-t+1)(s+2t+1)+\frac{t(t-3)}{2}-\frac{s(s+1)}{2} &> k=s+t\\
  \Leftrightarrow s^2-3t^2+2st+s-3t+2 &> 0.
 \end{align*}
 Since $t\leq s$, we can write $s=t+a$, where $a\geq0$, and substitute:
 \begin{align*}
  (t+a)^2-3t^2+2(t+a)t+(t+a)-3t+2 = a^2+a+4at-2t+2.
 \end{align*}
 The right hand side is $\leq0$ when $a=0$, and $>0$ when $a\geq1$. Thus $\binom{3}{k}_{s}-\binom{3}{k-(s+1)}_{s}>k$ whenever $k<2s$ and $\binom{3}{k}_{s}-\binom{3}{k-(s+1)}_{s}\leq k$ when $k=2s$.
 
 \item $N\geq4$: Recall that $\binom{N}{k}_{s}-\binom{N}{k-1}_{s}=\binom{N-1}{k}_{s}-\binom{N-1}{k-(s+1)}_{s}$ (see proof of \Cref{strongunimodality}). We will proceed by induction on N, and assume that $\binom{N-1}{k}_{s}-\binom{N-1}{k-(s+1)}_{s}>k$ for $k\leq\lfloor\frac{sN}{2}\rfloor$. Then $\binom{N}{k}_{s}-\binom{N}{k-1}_{s}>k$ and hence $\binom{N}{k}_{s}-\binom{N}{k-(s+1)}_{s}>k$ for $k\leq\lfloor\frac{sN}{2}\rfloor$.\\
 For the base case $N=4$, we only have $\binom{N-1}{k}_{s}-\binom{N-1}{k-(s+1)}_{s}>k$ for $k<\lfloor\frac{sN}{2}\rfloor$. So we do the case $k=\lfloor\frac{sN}{2}\rfloor=2s$ separately: From the proof of part (b), we have $\binom{4}{2s}_{s}-\binom{4}{2s-1}_{s}=\binom{3}{2s}_{s}-\binom{3}{2s-(s+1)}_{s}=3s+1+\frac{s(s-3)}{2}-\frac{s(s+1)}{2}=s+1$. Hence $\binom{4}{2s}_{s}-\binom{4}{s-1}_{s}=[\binom{4}{2s}_{s}-\binom{4}{2s-1}_{s}]+[\binom{4}{2s-1}_{s}-\binom{4}{2s-2}_{s}]+\ldots+[\binom{4}{s}_{s}-\binom{4}{s-1}_{s}]\stackrel{\text{(b)}}{>}s+1+
 \displaystyle\sum_{i=s}^{2s-1}i=s+1+\frac{(2s-1)(2s)}{2}-\frac{s(s-1)}{2}=\frac{s(3s+1)}{2}+1>2s$ 
 for $s\geq1$.\\
 Now for $N\geq4$, let $k=\lfloor\frac{sN}{2}\rfloor+t$, where $1\leq t<\lceil\frac{s}{2}\rceil$. Then
 \begin{align*}
 \binom{N}{k}_{s}&-\binom{N}{k-(s+1)}_{s}\stackrel{\text{sym.}}{=}\binom{N}{sN-k}_{s}-\binom{N}{k-(s+1)}_{s}\\
 &=\binom{N}{\lceil\frac{sN}{2}\rceil-t}_{s}-\binom{N}{\lfloor\frac{sN}{2}\rfloor+t-(s+1)}_{s}\\
 &=\binom{N}{\lceil\frac{sN}{2}\rceil-t}_{s}-\binom{N}{\lceil\frac{sN}{2}\rceil-t-1}_{s}+\\
 &\binom{N}{\lceil\frac{sN}{2}\rceil-t-1}_{s}-\binom{N}{\lceil\frac{sN}{2}\rceil-t-2}_{s}+\ldots+\\
 &\binom{N}{\lfloor\frac{sN}{2}\rfloor+t-s}_{s}-\binom{N}{\lfloor\frac{sN}{2}\rfloor+t-s-1}_{s}\\
 &>\displaystyle\sum_{i=\lfloor\frac{sN}{2}\rfloor+t-s}^{\lceil\frac{sN}{2}\rceil-t}i \text{ (by induction hypothesis)}\\
 &=\frac{(\lceil\frac{sN}{2}\rceil-t)(\lceil\frac{sN}{2}\rceil-t+1)}{2}-\frac{(\lfloor\frac{sN}{2}\rfloor+t-s)(\lfloor\frac{sN}{2}\rfloor+t-s-1)}{2}\\
 &=\frac{\lceil\frac{sN}{2}\rceil^2-\lfloor\frac{sN}{2}\rfloor^2+2s\lfloor\frac{sN}{2}\rfloor+sN-s^2-s-2tsN+2ts}{2}\\
 &=\begin{cases} (s^2N+sN-s^2-s-2tsN+2ts)/2 &\mbox{if } 2\mid sN \\
(s^2N+2sN-s^2-2s-2tsN+2ts)/2 & \mbox{if } 2\nmid sN \end{cases}\\
 &\stackrel{\text{?}}{>} \lfloor\frac{sN}{2}\rfloor+t=k.
 \end{align*}
 \underline{Case 1:} $2\mid sN$. Assume $2t\leq s-1$.
 \begin{align*}
 \frac{s^2N+sN-s^2-s-2tsN+2ts}{2}>\frac{sN}{2}+t\\
 \Leftrightarrow s^2N-s^2-s>2tsN-2ts+2t=2t(sN-s+1).
 \end{align*}
 For $N\geq4$, $2t\leq s-1$ implies
 \begin{align*}
 2t(sN-s+1)\leq s^2N-s^2-sN+2s-1<s^2N-s^2-s.
 \end{align*}
 \underline{Case 2:} $2\nmid sN$. Assume $2t\leq s$.
 \begin{align*}
 \frac{s^2N+2sN-s^2-2s-2tsN+2ts}{2}>\frac{sN-1}{2}+t\\
 \Leftrightarrow s^2N+sN-s^2-2s+1>2tsN-2ts+2t=2t(sN-s+1).
 \end{align*}
 For $N\geq4$, $2t\leq s$ implies
 \begin{align*}
 2t(sN-s+1)\leq s^2N-s^2+s<s^2N+sN-s^2-2s+1.
 \end{align*}
 So we have shown $\binom{N}{k}_{s}-\binom{N}{k-(s+1)}_{s}>k$ for\\
 $k\leq \begin{cases}
  \lfloor\frac{s(N+1)}{2}\rfloor-1 &\mbox{if } 2\mid s\\
  \lfloor\frac{s(N+1)}{2}\rfloor &\mbox{if } 2\nmid s \text{ and } 2\mid N\\
  \lfloor\frac{s(N+1)}{2}\rfloor-1 &\mbox{if } 2\nmid s \text{ and } 2\nmid N.
 \end{cases}$\\
 So we still need to show the case $k=\frac{s(N+1)}{2}$ if $2\mid s$ and if $2\nmid sN$. Assuming that
 $\binom{N-1}{k}_{s}-\binom{N-1}{k-1}_{s}=\binom{N-2}{k}_{s}-\binom{N-2}{k-(s+1)}_{s}>k$ for $k\leq\frac{s(N-1)}{2}$ (strong induction), we have
 \begin{align*}
 \binom{N}{\frac{s(N+1)}{2}}_{s}&-\binom{N}{\frac{s(N+1)}{2}-(s+1)}_{s}\stackrel{\text{sym.}}{=}
 \binom{N}{\frac{s(N-1)}{2}}_{s}-\binom{N}{\frac{s(N-1)}{2}-1}_{s}\\
 &=\binom{N-1}{\frac{s(N-1)}{2}}_{s}-\binom{N-1}{\frac{s(N-1)}{2}-(s+1)}_{s}\\
 &=\binom{N-1}{\frac{s(N-1)}{2}}_{s}-\binom{N-1}{\frac{s(N-1)}{2}-1}_{s}+\ldots+\\
 &\binom{N-1}{\frac{s(N-1)}{2}-s}_{s}-\binom{N-1}{\frac{s(N-1)}{2}-s-1}_{s}\\
 &>\displaystyle\sum_{i=\frac{s(N-1)}{2}-s}^{\frac{s(N-1)}{2}}i=\frac{s^2N+sN-2s^2-2s}{2}\\
 &\stackrel{\text{?}}{>}\frac{s(N+1)}{2}=k.
 \end{align*}
 The last inequality holds for $N\geq4$ and $s\geq2$. The case $s=1$ corresponds to the usual binomial
 coefficients and can be shown directly for $N\geq4$.\\
 For the base case $N=4$, we only have the weaker assumptions from parts (a) and (b). But from the 
 proof of part (b) we have
 \begin{align*}
 \binom{4}{\frac{5s}{2}}_{s}&-\binom{4}{\frac{5s}{2}-(s+1)}_{s}
 \stackrel{\text{sym.}}{=}\binom{4}{\frac{3s}{2}}_{s}-\binom{4}{\frac{3s}{2}-1}_{s}\\
 &=\binom{3}{\frac{3s}{2}}_{s}-\binom{3}{\frac{3s}{2}-(s+1)}_{s}\\
 &=(\frac{s}{2}+1)(2s+1)+\frac{\frac{s}{2}(\frac{s}{2}-3)}{2}-\frac{s(s+1)}{2}
 >\frac{5s}{2} \text{ for all } s.
 \end{align*}
 Also, for the case $N=5$, we only have $\binom{4}{2s}_{s}-\binom{4}{2s-1}_{s}=s+1$. So we have $\binom{5}{3s}_{s}-\binom{5}{3s-(s+1)}_{s}=\ldots>s+1+\displaystyle\sum_{i=s}^{2s-1}i=\frac{s(3s+1)}{2}+1>3s$ for $s\geq2$. For $s=1$, $\binom{5}{3}-\binom{5}{1}=5>3$.\\
 Now if $k=\lfloor\frac{s(N+1)}{2}\rfloor+1$ then
 \begin{align*}
 \binom{N}{\lfloor\frac{s(N+1)}{2}\rfloor+1}_{s}&-
 \binom{N}{\lfloor\frac{s(N+1)}{2}\rfloor+1-(s+1)}_{s}\\
 &\stackrel{\text{sym.}}{=}\binom{N}{\lceil\frac{s(N-1)}{2}\rceil-1}_{s}-
 \binom{N}{\lfloor\frac{s(N-1)}{2}\rfloor}_{s}
 \leq0 \leq k.
 \end{align*}
 Thus $k=\lfloor\frac{s(N+1)}{2}\rfloor+1$ is the smallest $k$ such that $\binom{N}{k}_{s}-\binom{N}{k-(s+1)}_{s}\leq k$.
 \end{enumerate}
\end{proof}

\begin{Corollary} ($c=2$)
If $f_1,\ldots,f_{n+2}\in K[x_1,\ldots,x_n]$ with $deg(f_i)=d$ $\forall i=1,\ldots,n+2$, then the 
smallest $D$ such that $\chi(D)\leq D$  is 
\begin{enumerate}[label=(\alph*),ref=(\alph*)]
 \item $\geq 2(d-1)$ for $n=2$,
 \item $\geq \lfloor\frac{(d-1)(n+2)}{2}\rfloor+1$ for $n\geq3$.
\end{enumerate}
\end{Corollary}

Based on the $c=1$ and $c=2$ cases, it appears as though the optimal starting value for $D$ in the XL algorithm will be approximately $(d-1)(n+c)/c$.

\section{Comparison with Degree of Regularity}

Some authors (for example see \cite{YCY13} or \cite{XLF4})
have compared the performance of XL and Groebner bases algorithms such as $F_4$ and $F_5$.
The conclusions seem to be that  each one can perform better than the other,
depending on the parameters and their
relationships.

The key concept in the complexity analysis of algorithms for computing Groebner bases is the degree of regularity, which is usually analyzed with Hilbert series. Since this paper also uses Hilbert series to analyze the maximal degree $D$ in the XL algorithm, it is reasonable to think about a comparison between the degree of regularity and the degree $D$ studied in this paper.

 In order to compare these two concepts, we will focus on the generic case. This is also similar to what is done in complexity analysis of Groebner basis algorithms. 
 Furthermore, we will just deal with the $c=1$ case.
 We will assume all polynomials have the same degree in order to get an interesting comparison.

Thus, suppose we have $n+1$ generic polynomials $g_1,\ldots,g_{n+1}\in K[x_1,\ldots,x_n]$ 
with $\deg(g_i)=d$ for all $i$.
Then \Cref{lb} tells us that the smallest $D$ such that the XL algorithm terminates is at most 
\[
(d-1)(n+1)-1
\] and we conjectured in \Cref{conj32} that this is exact.

Note that by \cite[Remark 1]{XLF4}, the behaviour of the XL algorithm on the homogenization of a polynomial system $V$ is the same as on the system $V$ itself.

On the other hand, consider the degree of regularity.
As remarked in \cite{Diem04}, the homogenization of this system forms a regular sequence. By \cite{BFSY05}, the degree of regularity of this homogenized system is exactly the Macaulay bound: $$\sum_{i=1}^{n+1}(d_i-1)+1.$$
When $d_i=d$ for all $i$, this is 
\[
(n+1)(d-1)+1. 
\]
Thus, our value for $D^*$ is   smaller than the degree of regularity by 2.

\bigskip

There is some doubt as to whether the degree of regularity is a good measure for complexity when the system is not a regular sequence, see \cite{DS13} and \cite{Elisa}. On the other hand, our $D^*$ determines the number of columns of the matrix calculated during the XL algorithm (and indirectly also the number of rows). The main computational part of the XL algorithm is Gaussian elimination on this matrix, and therefore the complexity of the XL algorithm is directly related to the complexity of Gaussian elimination, which increases with the matrix size. Thus, an upper bound on $D^*$ is a useful measure for the complexity of XL.

\bigskip

{\bf References}

\bigskip

\FloatBarrier
\bibliographystyle{alpha}
\bibliography{mybibXL}

\end{document}